\input amstex
\let\myfrac=\frac
\let\frac=\myfrac
\input epsf




\loadeufm \loadmsam \loadmsbm
\message{symbol names}\UseAMSsymbols\message{,}

\font\myfontdefault=cmr10

\font\mytdmchapfont=cmb10 at 14pt
\font\mytdmheadfont=cmb10 at 10pt
\font\mytdmsubheadfont=cmr10

\magnification 1200
\newif\ifinappendices
\newif\ifundefinedreferences
\newif\ifchangedreferences
\newif\ifloadreferences
\newif\ifmakebiblio
\newif\ifmaketdm

\undefinedreferencestrue
\changedreferencesfalse


\loadreferencestrue
\makebibliofalse
\maketdmfalse

\def\headpenalty{-400}     
\def\proclaimpenalty{-200} 

%
%

\def\alphanum#1{\ifcase #1 _\or A\or B\or C\or D\or E\or F\or G\or H\or I\or J\or K\or L\or M\or N\or O\or P\or Q\or R\or S\or T\or U\or V\or W\or X\or Y\or Z\fi}
\def\gobbleeight#1#2#3#4#5#6#7#8{}

\newwrite\references
\newwrite\tdm
\newwrite\biblio

\newcount\chapno
\newcount\headno
\newcount\subheadno
\newcount\procno
\newcount\figno
\newcount\citationno

\def\setcatcodes{%
\catcode`\!=0 \catcode`\\=11}%

\ifloadreferences
    {\catcode`\@=11 \catcode`\_=11%
    \global\def\_@citation@AndBarbBegZegh{1}
\global\def\_@citation@Ghomi{2}
\global\def\_@citation@GuanSpruck{3}
\global\def\_@citation@LabA{4}
\global\def\_@citation@Petersen{5}
\global\def\_@citation@Ros{6}
\global\def\_@citation@SmiSLC{7}
\global\def\_@citation@SmiFCS{8}
\global\def\_@citation@SmiNLD{9}
\global\def\_@citation@SmiCGC{10}
\global\def\_@citation@SmiPPH{11}
\global\def\_@citation@TrudWang{12}
\global\def\_@proc@ThmExistence{1.1}
\global\def\_@head@SpecialLagrangianCurvature{2}
\global\def\_@head@HeadConvexImmersions{3}
\global\def\_@proc@LemmaClosednessOfGraphProperty{3.1}
\global\def\_@proc@PropBasicResult{4.1}
\global\def\_@proc@PropConvexityOfDistanceFunction{4.2}
\global\def\_@proc@CorBoundedAwayFromSphere{4.3}
\global\def\_@proc@PropNonTrivialInterior{4.4}
\global\def\_@proc@LemmaPrecompactnessOfStrictlyConvexImmersions{5.1}
\global\def\_@proc@PropBoundsObtainedFromTheBoundary{5.2}
\global\def\_@proc@LemmaCompactness{6.1}
\global\def\_@proc@PropBasicExistenceResult{7.1}
    }%
\else
    \openout\references=references.tex
\fi

\newcount\newchapflag 
\newcount\showpagenumflag 

\global\chapno = -1 
\global\citationno=0
\global\headno = 0
\global\subheadno = 0
\global\procno = 0
\global\figno = 0

\def\resetcounters{%
\global\headno = 0%
\global\subheadno = 0%
\global\procno = 0%
\global\figno = 0%
}

\global\newchapflag=0 
\global\showpagenumflag=0 

\def\chinfo{\ifinappendices\alphanum\chapno\else\the\chapno\fi}%
\def\headinfo{\ifinappendices\alphanum\headno\else\the\headno\fi}%
\def\subheadinfo{\headinfo.\the\subheadno}
\def\procinfo{\headinfo.\the\procno}
\def\figinfo{\the\figno}        
\def\citationinfo{\the\citationno}%
\def\nextheadno{\global\advance\headno by 1 \global\subheadno = 0 \global\procno = 0}
\def\nextsubheadno{\global\advance\subheadno by 1}
\def\nextprocno{\global\advance\procno by 1 \procinfo}
\def\nextfigno{\global\advance\figno by 1 \figinfo}

{\global\let\noe=\noexpand%
%
%
\catcode`\@=11%
\catcode`\_=11%
\setcatcodes%
!global!def!_@@internal@@makeref#1{%
!global!expandafter!def!csname #1ref!endcsname##1{%
!csname _@#1@##1!endcsname%
!expandafter!ifx!csname _@#1@##1!endcsname!relax%
    !write16{#1 ##1 not defined - run saving references}%
    !undefinedreferencestrue%
!fi}}%
!global!def!_@@internal@@makelabel#1{%
!global!expandafter!def!csname #1label!endcsname##1{%
!edef!temptoken{!csname #1info!endcsname}%
!ifloadreferences%
    !expandafter!ifx!csname _@#1@##1!endcsname!relax%
        !write16{#1 ##1 not hitherto defined - rerun saving references}%
        !changedreferencestrue%
    !else%
        !expandafter!ifx!csname _@#1@##1!endcsname!temptoken%
        !else
            !write16{#1 ##1 reference has changed - rerun saving references}%
            !changedreferencestrue%
        !fi%
    !fi%
!else%
    !expandafter!edef!csname _@#1@##1!endcsname{!temptoken}%
    !edef!textoutput{!write!references{\global\def\_@#1@##1{!temptoken}}}%
    !textoutput%
!fi}}%
!global!def!makecounter#1{!_@@internal@@makelabel{#1}!_@@internal@@makeref{#1}}%
!unsetcatcodes%
}
\makecounter{ch}%
\makecounter{head}%
\makecounter{subhead}%
\makecounter{proc}%
\makecounter{fig}%
\makecounter{citation}%
\def\newref#1#2{%
\def\temptext{#2}%
\edef\bibliotextoutput{\expandafter\gobbleeight\meaning\temptext}%
\global\advance\citationno by 1\citationlabel{#1}%
\ifmakebiblio%
    \edef\fileoutput{\write\biblio{\noindent\hbox to 0pt{\hss$[\the\citationno]$}\hskip 0.2em\bibliotextoutput\medskip}}%
    \fileoutput%
\fi}%
\def\cite#1{%
$[\citationref{#1}]$%
\ifmakebiblio%
    \edef\fileoutput{\write\biblio{#1}}%
    \fileoutput%
\fi%
}%
%
%
%

\let\mypar=\par


\def\raggedleft{\leftskip=0pt plus 1fil \parfillskip=0pt}


\font\lettrinefont=cmr10 at 28pt
\def\lettrine #1[#2][#3]#4%
{\hangafter -#1 \hangindent #2
\noindent\hskip -#2 \vtop to 0pt{
\kern #3 \hbox to #2 {\lettrinefont #4\hss}\vss}}

\font\mylettrinefont=cmr10 at 28pt
\def\mylettrine #1[#2][#3][#4]#5%
{\hangafter -#1 \hangindent #2
\noindent\hskip -#2 \vtop to 0pt{
\kern #3 \hbox to #2 {\mylettrinefont #5\hss}\vss}}


\edef\Pagetitle={Blank}

\headline={\hfil\Pagetitle\hfil}

\footline={\hfil\myfontdefault\folio\hfil}

\def\nextoddpage
{
\newpage%
\ifodd\pageno%
\else%
    \global\showpagenumflag = 0%
    \null%
    \vfil%
    \eject%
    \global\showpagenumflag = 1%
\fi%
}


\def\newchap#1#2%
{%
%
%
\global\advance\chapno by 1%
\resetcounters%
%
%
\newpage%
\ifodd\pageno%
\else%
    \global\showpagenumflag = 0%
    \null%
    \vfil%
    \eject%
    \global\showpagenumflag = 1%
\fi%
\global\newchapflag = 1%
\global\showpagenumflag = 1%
%
%
{\font\chapfontA=cmsl10 at 30pt%
\font\chapfontB=cmsl10 at 25pt%
\null\vskip 5cm%
{\chapfontA\raggedleft\hfil%
{%
\ifnum\chapno=0
    \phantom{%
    \ifinappendices%
        Annexe \alphanum\chapno%
    \else%
        \the\chapno%
    \fi}%
\else%
    \ifinappendices%
        Annexe \alphanum\chapno%
    \else%
        \the\chapno%
    \fi%
\fi%
}%
\par}%
\vskip 2cm%
{\chapfontB\raggedleft%
\lineskiplimit=0pt%
\lineskip=0.8ex%
\hfil #1\par}%
\vskip 2cm%
}%
\edef\Pagetitle{#2}%
%
%
\ifmaketdm%
    \def\temp{#2}%
    \def\tempbis{\nobreak}%
    \edef\chaptitle{\expandafter\gobbleeight\meaning\temp}%
    \edef\mynobreak{\expandafter\gobbleeight\meaning\tempbis}%
    \edef\textoutput{\write\tdm{\bigskip{\noexpand\mytdmchapfont\noindent\chinfo\ - \chaptitle\hfill\noexpand\folio}\par\mynobreak}}%
\fi%
\textoutput%
}


\def\newhead#1%
{%
\ifhmode%
    \mypar%
\fi%
\ifnum\headno=0%
\ifinappendices
    \nobreak\vskip -\lastskip%
    \nobreak\vskip .5cm%
\fi
\else%
    \nobreak\vskip -\lastskip%
    \nobreak\vskip .5cm%
\fi%
\nextheadno%
\ifmaketdm%
    \def\temp{#1}%
    \edef\sectiontitle{\expandafter\gobbleeight\meaning\temp}%
    \edef\textoutput{\write\tdm{\noindent{\noexpand\mytdmheadfont\quad\headinfo\ - \sectiontitle\hfill\noexpand\folio}\par}}%
    \textoutput%
\fi%
\font\headfontA=cmbx10 at 14pt%
{\headfontA\noindent\headinfo\ - #1.\hfil}%
\nobreak\vskip .5cm%
}%


\def\newsubhead#1%
{%
\ifhmode%
    \mypar%
\fi%
\ifnum\subheadno=0%
\else%
    \penalty\headpenalty\vskip .4cm%
\fi%
\nextsubheadno%
\ifmaketdm%
    \def\temp{#1}%
    \edef\subsectiontitle{\expandafter\gobbleeight\meaning\temp}%
    \edef\textoutput{\write\tdm{\noindent{\noexpand\mytdmsubheadfont\quad\quad\subheadinfo\ - \subsectiontitle\hfill\noexpand\folio}\par}}%
    \textoutput%
\fi%
\font\subheadfontA=cmsl10 at 12pt
{\subheadfontA\noindent\subheadinfo\ #1.\hfil}%
\nobreak\vskip .25cm %
}%

%
%


\font\mathromanten=cmr10
\font\mathromanseven=cmr7
\font\mathromanfive=cmr5
\newfam\mathromanfam
\textfont\mathromanfam=\mathromanten
\scriptfont\mathromanfam=\mathromanseven
\scriptscriptfont\mathromanfam=\mathromanfive
\def\mathroman{\fam\mathromanfam}


\font\sf=cmss12

\font\sansseriften=cmss10
\font\sansserifseven=cmss7
\font\sansseriffive=cmss5
\newfam\sansseriffam
\textfont\sansseriffam=\sansseriften
\scriptfont\sansseriffam=\sansserifseven
\scriptscriptfont\sansseriffam=\sansseriffive
\def\mathsf{\fam\sansseriffam}


\font\bftwelve=cmb12

\font\boldten=cmb10
\font\boldseven=cmb7
\font\boldfive=cmb5
\newfam\mathboldfam
\textfont\mathboldfam=\boldten
\scriptfont\mathboldfam=\boldseven
\scriptscriptfont\mathboldfam=\boldfive
\def\mathbf{\fam\mathboldfam}


\font\mycmmiten=cmmi10
\font\mycmmiseven=cmmi7
\font\mycmmifive=cmmi5
\newfam\mycmmifam
\textfont\mycmmifam=\mycmmiten
\scriptfont\mycmmifam=\mycmmiseven
\scriptscriptfont\mycmmifam=\mycmmifive

\def\hexa#1{\ifcase #1 0\or 1\or 2\or 3\or 4\or 5\or 6\or 7\or 8\or 9\or A\or B\or C\or D\or E\or F\fi}
\mathchardef\mathi="7\hexa\mycmmifam7B
\mathchardef\mathj="7\hexa\mycmmifam7C


\font\mymsbmten=msbm10 at 8pt
\font\mymsbmseven=msbm7 at 5.6pt
\font\mymsbmfive=msbm5 at 4pt
\newfam\mymsbmfam
\textfont\mymsbmfam=\mymsbmten
\scriptfont\mymsbmfam=\mymsbmseven
\scriptscriptfont\mymsbmfam=\mymsbmfive

\mathchardef\mybeth="7\hexa\mymsbmfam69
\mathchardef\mygimmel="7\hexa\mymsbmfam6A
\mathchardef\mydaleth="7\hexa\mymsbmfam6B


\def\placelabel[#1][#2]#3{{%
\setbox10=\hbox{\raise #2cm \hbox{\hskip #1cm #3}}%
\ht10=0pt%
\dp10=0pt%
\wd10=0pt%
\box10}}%


\newif\ifinproclaim%
\global\inproclaimfalse%
\def\proclaim#1{%
\medskip%
%
%
\bgroup%
\inproclaimtrue%
\setbox10=\vbox\bgroup\leftskip=0.8em\noindent{\bftwelve #1}\sf%
}

\def\endproclaim{%
\egroup%
\setbox11=\vtop{\noindent\vrule height \ht10 depth \dp10 width 0.1em}%
\wd11=0pt%
\setbox12=\hbox{\copy11\kern 0.3em\copy11\kern 0.3em}%
\wd12=0pt%
\setbox13=\hbox{\noindent\box12\box10}%
\noindent\unhbox13%
\egroup%
\medskip\ignorespaces%
}

\def\proclaim#1{%
\medskip%
\bgroup%
\inproclaimtrue%
\noindent{\bftwelve #1}%
\nobreak\medskip%
\sf%
}

\def\endproclaim{%
\mypar\egroup\penalty\proclaimpenalty\medskip\ignorespaces%
}

\def\noskipproclaim#1{%
\medskip%
\bgroup%
\inproclaimtrue%
\noindent{\bf #1}\nobreak\sl%
}

\def\endnoskipproclaim{%
\mypar\egroup\penalty\proclaimpenalty\medskip\ignorespaces%
}


\def\ninn{{n\in\Bbb{N}}}

\def\proof{{\noindent\bf Proof:\ }}

\def\remark{{\noindent\sl Remark:\ }}
\def\example{{\noindent\sl Example:\ }}

\def\minf{\mathop{{\mathroman Inf}}}
\def\msf#1{{\mathsf #1}}

\def\qed{~$\square$}
\def\munion{\mathop{\cup}}
\def\minter{\mathop{\cap}}
\def\myitem#1{%
    \noindent\hbox to .5cm{\hfill#1\hss}
}

\catcode`\@=11
\def\Eqalign#1{\null\,\vcenter{\openup\jot\m@th\ialign{%
\strut\hfil$\displaystyle{##}$&$\displaystyle{{}##}$\hfil%
&&\quad\strut\hfil$\displaystyle{##}$&$\displaystyle{{}##}$%
\hfil\crcr #1\crcr}}\,}
\catcode`\@=12

\def\makeop#1{%
\global\expandafter\def\csname op#1\endcsname{{\mathroman #1}}}%

\def\makeopsmall#1{%
\global\expandafter\def\csname op#1\endcsname{{\mathroman{\lowercase{#1}}}}}%

\makeopsmall{ArcTan}%
\makeopsmall{ArcCos}%
\makeop{Arg}%
\makeop{Det}%
\makeop{Log}%
\makeop{Re}%
\makeop{Im}%
\makeop{Dim}%
\makeopsmall{Tan}%
\makeop{Ker}%
\makeopsmall{Cos}%
\makeopsmall{Sin}%
\makeop{Exp}%
\makeopsmall{Tanh}%
\makeop{Tr}%
\makeop{End}%
\makeop{Long}%
\makeop{Ch}%
\makeop{Exp}%
\makeop{Eval}%
\makeop{Lift}%
\makeop{Int}%
\makeop{Ext}%
\makeop{Aire}%
\makeop{Im}%
\makeop{Conf}%
\makeop{Exp}%
\makeop{Mod}%
\makeop{Log}%
\makeop{Sgn}%
\makeop{Ext}%
\makeop{Int}%
\makeop{Ln}%
\makeop{Dist}%
\makeop{Aut}%
\makeop{Id}%
\makeop{GL}%
\makeop{SO}%
\makeop{Homeo}%
\makeop{Vol}%
\makeop{Ric}%
\makeop{Hess}%
\makeop{Euc}%
\makeop{Isom}%
\makeop{Max}%
\makeop{SW}%
\makeop{SL}%
\makeop{Long}%
\makeop{Fix}%
\makeop{Wind}%
\makeop{Diag}%
\makeop{dVol}%
\makeop{Supp}%
\makeop{Symm}%
\makeop{Ad}%
\makeop{Diam}%
\makeop{loc}%
\makeopsmall{Sinh}%
\makeopsmall{Cosh}%
\makeop{Len}%
\makeop{Length}%
\makeop{Conv}%
\makeop{Min}%
\makeop{Area}%

\font\mycirclefont=cmsy7
\def\textcircle{{\raise 0.3ex \hbox{\mycirclefont\char'015}}}

\let\emph=\bf

\hyphenation{quasi-con-formal}

%
%

\ifmakebiblio%
    \openout\biblio=biblio.tex %
    {%
        \edef\fileoutput{\write\biblio{\bgroup\leftskip=2em}}%
        \fileoutput
    }%
\fi%

\newref{AndBarbBegZegh}{Andersson L., Barbot T., B\'eguin F., Zeghib A., Cosmological time versus CMC time in spacetimes of constant curvature, arXiv:math/0701452}
\newref{Ghomi}{Alexander S., Ghomi M., Wong J., Topology of Riemannian submanifolds with prescribed boundary, {\sl to appear in Duke Math. J.}}
\newref{GuanSpruck}{Guan B., Spruck J., The existence of hypersurfaces of constant Gauss curvature with prescribed boundary, {\sl J. Differential Geom.} {\bf 62} (2002), no. 2, 259--287}
\newref{LabA}{Labourie F., Un lemme de Morse pour les surfaces convexes (French), {\sl Invent. Math.} {\bf 141} (2000), no. 2, 239--297}
\newref{Petersen}{Petersen P., {\sl Riemannian geometry}, Graduate Texts in Mathematics, {\bf 171}, Springer, New York, (2006)}
\newref{Ros}{Rosenberg H., Hypersurfaces of constant curvature in space forms, {\sl Bull. Sci. Math.} {\bf 117} (1993), no. 2, 211--239} 
\newref{SmiSLC}{Smith G., Special Lagrangian curvature, arXiv:math/0506230}
\newref{SmiFCS}{Smith G., Moduli of Flat Conformal Structures of Hyperbolic Type, arXiv:0804.0744}
\newref{SmiNLD}{Smith G., The Non-Linear Dirichlet Problem in Hadamard Manifolds, arXiv:0908.3590}
\newref{SmiCGC}{Smith G., Constant Gaussian Curvature Hypersurfaces in Hadamard Manifold, arXiv:0912.0248}
\newref{SmiPPH}{Smith G., The Plateau Problem in Hadamard Manifolds, arXiv:1002.2982}
\newref{TrudWang}{Trudinger N. S., Wang X., On locally locally convex hypersurfaces with boundary, {\sl J. Reine Angew. Math.} {\bf 551} (2002), 11--32}

\ifmakebiblio%
    {\edef\fileoutput{\write\biblio{\egroup}}%
    \fileoutput}%
\fi%

%
%
%
\document
\myfontdefault
\global\chapno=1
\global\showpagenumflag=1
\def\Pagetitle{}
\null
\vfill
\def\centre{\rightskip=0pt plus 1fil \leftskip=0pt plus 1fil \spaceskip=.3333em \xspaceskip=.5em \parfillskip=0em \parindent=0em}%
\def\textmonth#1{\ifcase#1\or January\or Febuary\or March\or April\or May\or June\or July\or August\or September\or October\or November\or December\fi}
\font\abstracttitlefont=cmr10 at 14pt
{\abstracttitlefont\centre The non-linear Plateau problem in non-positively curved manifolds\par}
\bigskip
{\centre Graham Smith\par}
\bigskip
{\centre 2nd April 2010\par}
\bigskip
{\centre IMPA,\par
Estrada Dona Castorina 110,\par
Rio de Janeiro,\par
BRASIL 22460-320\par}
\bigskip
\noindent{\emph Abstract:\ }Using the Perron method, we prove the existence of hypersurfaces of prescribed special Lagrangian curvature with prescribed boundary inside complete Riemannian manifolds of non-positive curvature.
%
%
\par 
\vfill
\nextoddpage
\global\pageno=1
\myfontdefault
\def\Pagetitle{\sl The non-linear Plateau problem in non-positively curved manifolds}
\newhead{Introduction}
\noindent In this paper we prove the existence of solutions of the Plateau problem for hypersurfaces of prescribed curvature and prescribed boundary in manifolds of non-positive sectional curvature. The curvature notion used here - special Lagrangian curvature - was introduced by the author in \cite{SmiSLC} and constitutes a higher dimensional generalisation of two dimensional extrinsic curvature. Its interest is two-fold: firstly, it possesses strong regularity properties that translate into very simple geometric behaviour of limits of sequences of hypersurfaces of prescribed special Lagrangian curvature; and, secondly, in low dimensions, it reduces to certain well known notions of curvature. Explicitely, let $M^{n+1}$ be a Riemannian manifold and let $N^n$ be a locally strictly convex, smooth, immersed hypersurface in $M$. The special Lagrangian curvature, $R_\theta(N)$, is a function of the second fundamental form of $N$ which depends on an angle parameter, $\theta\in[0,n\pi/2[$ (see Section \headref{SpecialLagrangianCurvature} for more details). In low dimensions, we have:
\medskip
\myitem{(i)} when the ambient manifold is $3$-dimensional, and thus when $N$ is a surface:
$$
R_{\pi/2} = K_e^{1/2},
$$
\noindent where $K_e$ is the extrinsic curvature of $N$; and
\medskip
\myitem{(ii)} when the ambient manifold is $4$-dimensional:
$$
R_\pi = (K_e/H)^{1/2},
$$
\noindent where $H$ is the mean curvature of $N$. Moreover:
$$
R_{\pi/2} = (R^N - R^M + \opRic^M(\msf{N},\msf{N}))^{1/2},
$$
\noindent where $R^N$ and $R^M$ are the scalar curvatures of $N$ and $M$ respectively, $\opRic^M$ is the Ricci curvature of $M$ and $\msf{N}$ is the unit exterior normal vector over $N$.
\medskip
\noindent Any existence result for hypersurfaces of prescribed special Lagrangian curvature thus translates in particular into existence results for the three notions of curvature given above. 
\proclaim{Theorem \nextprocno}
\noindent Let $M^{n+1}$ be a complete $(n+1)$-dimensional Riemannian manifold of non-positive sectional curvature. Let $\phi:M\rightarrow]0,\infty[$ be a strictly positive smooth function. Let $\theta\in[(n-1)\pi/2,n\pi/2[$ be an angle, and let $(N^n,\partial N^n)$ be a compact locally convex immersed hypersurface in $M$ such that, for all $p\in N$:
$$
R_\theta(N)(p) \geqslant \phi(p),
$$
\noindent in the weak sense. If $\theta>(n-1)\pi/2$, then there exists a compact, locally convex immersed hypersurface, $(\hat{N}^n,\partial\hat{N}_n)$ in $M$ such that:
\medskip
\myitem{(i)} $\hat{N}$ has smooth interior and is $C^{0,1}$ up to its boundary;
\medskip
\myitem{(ii)} $\partial\hat{N}=\partial N$ and $\hat{N}$ is bounded by $N$; and
\medskip
\myitem{(iii)} for all $p\in\hat{N}$:
$$
R_\theta(\hat{N})(p) = \phi(p).
$$
\noindent In fact, $(\hat{N},\partial\hat{N})$ is isotopic by locally convex immersions to $(N,\partial N)$.
\medskip
\noindent If $\theta=(n-1)\pi/2$, then the same result holds provided that, in addition, the shape operator of $N$ is everywhere bounded below by $\epsilon\opId$ in the weak sense, for some $\epsilon>0$.
\endproclaim
\proclabel{ThmExistence}
\remark See Sections \headref{SpecialLagrangianCurvature} and \headref{HeadConvexImmersions} for terminology.
\medskip
\remark This generalises Theorem $1.1$ of \cite{SmiNLD} to more general manifolds and also to the case of submanifolds with boundary.
\medskip
\remark When the ambient manifold is $3$-dimensional, and when $N$ is thus a surface, Theorem \procref{ThmExistence} constitutes a generalisation of the existence part of Proposition $5.0.3$ of \cite{LabA}, which itself constitutes the analytic core of that paper (see also \cite{SmiPPH}).
\medskip
\remark When the ambient manifold is $4$-dimensional, although this result proves existence of hypersurfaces for the case where $R_\pi=(K_e/H)^{1/2}$ is prescribed, it says nothing about the case of $R_{\pi/2}$. However, when $M$ is hyperbolic, we can canonically associate to it a de-Sitter manifold, $M'$, which is dual to $M$ in a certain sense. This duality has the effect of interchanging $R_\pi$ and $R_{\pi/2}$. Moreover, since $R^{M'}$ and $\opRic^{M'}$ are constant in this case, this yields an existence result for $3$-dimensional hypersurfaces of prescribed scalar curvature and prescribed boundary in $4$-dimensional de-sitter Manifolds (see also \cite{AndBarbBegZegh}).
\medskip
\remark In general, the Plateau problem is stated as follows: given a Riemannian manifold, $M$, and a finite family $\Gamma_1,...,\Gamma_n$ of smooth, codimension $2$ submanifolds in $M$, when does there exists an immersed hypersurface $N\subseteq M$ of constant curvature spanning $\Gamma_1\munion...\munion\Gamma_n$, i.e. such that:
$$
\partial N = \Gamma_1\munion...\munion\Gamma_n?
$$
\noindent This result therefore reduces the Plateau problem for prescribed special Lagrangian curvature to the weaker problem of determining when there exists a locally convex immersed hypersurface of $M$ spanning $\Gamma_1\munion...\munion\Gamma_n$. This problem is addressed, for example, by Alexander, Ghomi and Wong in \cite{Ghomi} and Rosenberg in \cite{Ros}.
\medskip
\noindent The proof of this result uses the Perron method, first applied to the study of hypersurfaces in two independent and simultaneous papers by Gaun and Spruck, \cite{GuanSpruck}, and Trudinger and Wang, \cite{TrudWang}, where they prove the existence of hypersurfaces of constant extrinsic curvature and prescribed boundary in $\Bbb{R}^{n+1}$. The main new ingredients used here involve various recent results by the author - \cite{SmiSLC}, \cite{SmiNLD}, \cite{SmiCGC} and \cite{SmiPPH} - of which the most significant are the regularity properties for special Lagrangian curvature (see \cite{SmiSLC}), and a regularity result near the boundary for locally convex immersions in Riemannian manifolds (see \cite{SmiPPH}). These results are here combined around the following compactness result for locally strictly convex immersions in Riemannian manifolds, which generalises the Main Lemma of \cite{TrudWang} and is of independant interest:
\goodbreak
\proclaim{Lemma \procref{LemmaCompactness}}
\noindent Let $M^{n+1}$ be an $(n+1)$-dimensional Riemannian manifold. Choose $\epsilon>0$ and let $(\Sigma_n,\partial\Sigma_n)_\ninn$ be a sequence of compact, locally convex immersed hypersurfaces such that:
\medskip
\myitem{(i)} $\Gamma_n:=\partial\Sigma_n$ is $C^\infty$; and
\medskip
\myitem{(ii)} the shape operator of $\Sigma_n$ is greater than $\epsilon\opId$ in the weak sense.
\medskip
\noindent Let $\Gamma_0$ be a smooth, compact, codimension $2$ submanifold of $M$ and suppose that $(\Gamma_n)_\ninn$ converges to $\Gamma_0$ in the $C^\infty$ sense. If there exists a compact subset $K\subseteq M$ and a real number $B>0$ such that, for all $n$:
$$
\Sigma_n\subseteq K,\qquad \opVol(\Sigma_n)\leqslant B,
$$
\noindent then there exists a $C^{0,1}$ locally convex immersed hypersurface $\Sigma_0$ in $M$ such that:
\medskip
\myitem{(i)} $\Gamma_0=\partial\Sigma_0$;
\medskip
\myitem{(ii)} the shape operator of $\Sigma_0$ is greater than $\epsilon\opId$ in the weak sense; and
\medskip
\myitem{(iii)} $(\Sigma_n)_\ninn$ subconverges to $\Sigma_0$.
\endproclaim
\noindent This paper is structured as follows: Sections $2$ and $3$ contain definitions and notation; Sections $4$, $5$ and $6$ are devoted to the proof of Lemma \procref{LemmaCompactness}; Section $7$ provides a slight generalisation of the solution to the Dirichlet problem studied in \cite{SmiNLD}, which forms the analytic core of the Perron Method in this paper; and in Section $8$ we prove Theorem \procref{ThmExistence}.
\medskip
\noindent The author would like to thank the Instito Nacional de Matem\'atica Pura e Aplicada in Rio de Janeiro, Brazil for providing the conditions required to complete this paper.
\newhead{Immersed Submanifolds and Special Lagrangian Curvature}
\noindent Let $M^{n+1}$ be an $(n+1)$-dimensional Riemannian manifold. An {\bf immersed submanifold} is a pair $\Sigma=(S,i)$ where $S$ is a smooth manifold and $i:S\rightarrow M$ is a smooth immersion. An {\bf immersed hypersurface} is an immersed submanifold of codimension $1$. We say that an immersed hypersurface is locally (strictly) convex if and only if its shape operator is everywhere positive definite. The special Lagrangian curvature, which is only defined for locally strictly convex immersed hypersurfaces, is defined as follows (see \cite{SmiSLC} for details): denote by $\opSymm(\Bbb{R}^n)$ the space of symmetric matrices over $\Bbb{R}^n$. We define $\Phi:\opSymm(\Bbb{R}^n)\rightarrow\Bbb{C}^*$ by:
\headlabel{SpecialLagrangianCurvature}
$$
\Phi(A) = \opDet(I+iA).
$$
\noindent Since $\Phi$ never vanishes and $\opSymm(\Bbb{R}^n)$ is simply connected, there exists a unique analytic function $\tilde{\Phi}:\opSymm(\Bbb{R}^n)\rightarrow\Bbb{C}$
such that:
$$
\tilde{\Phi}(I) = 0,\qquad e^{\tilde{\Phi}(A)} = \Phi(A)\qquad\forall A\in\opSymm(\Bbb{R}^n).
$$
\noindent We define the function $\opArcTan:\opSymm(\Bbb{R}^n)\rightarrow(-n\pi/2,n\pi/2)$ by:
$$
\arctan(A) = \opIm(\tilde{\Phi}(A)).
$$
\noindent This function is trivially invariant under the action of $O(\Bbb{R}^n)$. Moreover, if $\lambda_1, ..., \lambda_n$ are the eigenvalues of $A$, then:
$$
\opArcTan(A) = \sum_{i=1}^n\opArcTan(\lambda_i).
$$
\noindent For $r>0$, we define:
$$
\opSL_r(A) = \opArcTan(r^{-1}A).
$$
\noindent If $A$ is positive definite, then $SL_r$ is a strictly decreasing function of $r$. Moreover, $SL_\infty=0$ and $SL_0=n\pi/2$. Thus, for all $\theta\in]0,n\pi/2[$, there exists a unique $r>0$ such that:
$$
SL_r(A) = \theta.
$$
\noindent We define $R_\theta(A) = r$. $R_\theta$ is also invariant under the action of $O(n)$ on the space of positive definite, symmetric matrices.
\medskip
\noindent Let $M^{n+1}$ be an oriented Riemannian manifold of dimension $n+1$. Let $\Sigma=(S,i)$ be a locally strictly convex, immersed hypersurface in $M$. For $\theta\in]0,n\pi/2[$, we define $R_\theta(\Sigma)$ (the {\emph $\theta$-special Lagrangian curvature} of $\Sigma$) by:
$$
R_\theta(\Sigma) = R_\theta(A_\Sigma),
$$
\noindent where $A_\Sigma$ is the shape operator of $\Sigma$.
\goodbreak
\newhead{Locally Convex Hypersurfaces}
\noindent Let $M^{n+1}$ be a Riemannian manifold. A {\bf locally convex hypersurface} in $M$ is a pair $\Sigma=(i,S^n)$ where $S$ is an $n$-dimensional topological manifold and $i:S\rightarrow M$ is a continuous map such that, for all $p\in S$, there exists a neighbourhood, $U$, of $p$ in $S$, a convex subset $K\subseteq M$ with non-trivial interior, and an open subset $V\subseteq\partial K$ such that $i$ restricts to a homeomorphism from $U$ to $V$. We refer to such a triplet $(U,V,K)$ as a {\bf convex chart} of $\Sigma$. Pulling back the metric on $M$ through $i$ yields a natural length metric on $\Sigma$ which we denote by $d_\Sigma$. Let $(\Sigma_n)_\ninn = (i_n,S_n)_\ninn$ and $S_0 = (i_0,S_0)$ be convex immersions. We say that $(\Sigma_n)_\ninn$ {\bf converges} to $\Sigma_0$ if and only if:
\headlabel{HeadConvexImmersions}
\medskip
\myitem{(i)} $(S_n,d_{\Sigma_n})_\ninn$ converges to $(S_0,d_{\Sigma_0})$ in the Gromov-Hausdorff sense; and
\medskip
\myitem{(ii)} $(i_n)_\ninn$ converges to $i_0$ locally uniformly.
\medskip
\noindent Let $\Sigma=(i,S)$ and $\Sigma'=(i',S')$ be two locally convex hypersurfaces in $M$. We say that $\Sigma$ and $\Sigma'$ are equivalent if and only if there exists a homeomorphism $\phi:S\rightarrow S'$ such that:
$$
i = i'\circ\phi.
$$
\example Let $K\subseteq M$ be a convex subset with non trivial interior. Then any open subset of $\partial K$ is a locally convex hypersurface.\qed
\medskip
\example Let $\Sigma$ be a (smooth) hypersurface on $M$. $\Sigma$ is a locally convex hypersurface if and only if its second fundamental form is everywhere non-negative definite.\qed
\medskip
\noindent For $\epsilon>0$, we say that the second fundamental form of $\Sigma=(i,S)$ is {\bf bounded below} by $\epsilon\opId$ in the {\bf weak sense} if and only if, for every $p\in S$, and for every supporting normal, $\msf{N}_p$ of $S$ at $p$, there exists a smooth, embedded hypersurface $\Sigma'$ such that:
\medskip
\myitem{(i)} $\Sigma'$ passes through $p$;
\medskip
\myitem{(ii)} the normal to $\Sigma'$ at $p$ is $\msf{N}_p$;
\medskip
\myitem{(iii)} the shape operator of $\Sigma'$ at $p$ is equal to $\epsilon\opId$; and
\medskip
\myitem{(iii)} there exists a neighbourhood $U$ of $p$ in $S$ such that $i(U)$ lies entirely to one side of $\Sigma'$.
\medskip
\noindent Observe that this definition is stronger than the usual weak notion of curvature bounds (in the Alexandrov, or viscosity, sense), since the latter does not assume the existence of bounding hypersurfaces normal to arbitrary directions at arbitrary points.
\medskip
\noindent We may define the notion of weak lower (and upper) bounds for the special Lagrangian curvature in an analogous way.
\medskip
\noindent Suppose now that $M$ is a Hadamard manifold. Let $K\subseteq M$ be a convex set with non-trivial interior. Let $K^o$ be the interior of $K$. We define $\pi_K:M\setminus K^o\rightarrow \partial K$ to be projection onto the closest point in $\partial K$. Let $V\subseteq\partial K$. We call the set $\pi_K^{-1}(V)$ the {\bf end} of $V$, and we denote it by $\Cal{E}(V)$. Trivially, $\Cal{E}(V)$ is foliated by half geodesics leaving points in $V$ in directions normal to $K$. Let $\Sigma$ be a locally convex hypersurface. Let $(U,V,K)$ and $(U',V',K')$ be convex charts of $\Sigma$. Trivially:
$$
\pi_K^{-1}(i(U\minter U')) = \pi_{K'}^{-1}(i(U\minter U')).
$$
\noindent We thus define the {\bf end} of $\Sigma$ to be the manifold (with non-smooth, concave boundary) whose coordinate charts are the ends of the convex charts of $\Sigma$. We denote this manifold by $\Cal{E}(\Sigma)$. $\Cal{E}(\Sigma)$ has the following properties:
\medskip
\myitem{(i)} $\Sigma$ naturally embeds as the boundary of $\Cal{E}(\Sigma)$;
\medskip
\myitem{(ii)} in the complement of $\Sigma$, $\Cal{E}(\Sigma)$ has the structure of a smooth Riemannian manifold of non-positive curvature;
\medskip
\myitem{(iii)} $\Cal{E}(\Sigma)$ is foliated by half geodesics leaving points in $\Sigma$ in directions normal to $\Sigma$; and
\medskip
\myitem{(iv)} there exists a natural embedding $I:\Cal{E}(\Sigma)\rightarrow M$ which restricts to $i$ over $\Sigma$ and which is a local diffeomorphism over the complement of $\Sigma$.
\medskip
\noindent We say that a subset $K\subseteq\Cal{E}(\Sigma)$ is {\bf semi-convex} if and only if for every geodesic segment $\gamma:[0,1]\rightarrow\Cal{E}(\Sigma)$ contained within $\Cal{E}(\Sigma)$, if $\gamma(0),\gamma(1)\in K$, then the whole of $\gamma$ is contained in $K$. Let $K$ be a semi-convex subset of the end of $\Sigma$ which contains $\Sigma$ and coincides with $\Sigma$ outside a convex set. $(\partial K,I|_{\partial K})$ defines a convex immersion in $M$ which, by abuse of notation, we simply denote by $\partial K$ (see Section $3$ of \cite{SmiPPH} for details). Let $\Sigma$ and $\Sigma'$ be two locally convex hypersurfaces in $M$. We say that $\Sigma$ is {\bf bounded by} $\Sigma'$ (and $\Sigma'$ {\bf bounds} $\Sigma$) if and only if there exists a semi-convex subset, $K\subseteq\Cal{E}(\Sigma)$, which contains $\Sigma$ and which coincides with $\Sigma$ outside a compact set such that $\Sigma'$ is equivalent to $\partial K$. In this case, we often identify $\Sigma'$ with $\partial K$ and thus view it as a subset of $\Cal{E}(\Sigma)$.
\medskip
\example Let $K,K'\subseteq M$ be two convex sets. Then $\partial K$ is bounded by $\partial K'$ if and only if $K\subseteq K'$.\qed
\medskip
\noindent The property of boundedness is preserved by passage to limits. Indeed, if $(\Sigma,\partial\Sigma)$ is a locally convex immersion with boundary, then, for all $r>0$, we define $B_r(\partial\Sigma)$ to be the set of all points in $\Sigma$ whose intrinsic distance to $\partial\Sigma$ is less than $r$. We obtain:
\proclaim{Lemma \nextprocno}
\noindent Let $(\Sigma_n)_\ninn,\Sigma_0$ and $(\Sigma'_n)_\ninn,\Sigma'_0$ be compact, convex immersed hypersurfaces in $M$ such that $(\Sigma_n)_\ninn$ and $(\Sigma'_n)_\ninn$ converge to $\Sigma_0$ and $\Sigma'_0$ respectively. Suppose that:
\medskip
\myitem{(i)} for all $n>0$, $\Sigma'_n$ bounds $\Sigma_n$; and
\medskip
\myitem{(ii)} there exists $r>0$ such that, for all $n\in\Bbb{N}$:
$$
B_r(\partial\Sigma_n) = B_r(\partial\Sigma_n').
$$ 
\noindent Then $\Sigma'_0$ also bounds $\Sigma_0$.
\endproclaim
\proclabel{LemmaClosednessOfGraphProperty}
\proof See Lemma $3.2$ of \cite{SmiPPH}.\qed
\goodbreak
\newhead{Strictly Convex Hypersurfaces of Euclidean Space}
\noindent Let $B_1(0)$ be the ball of radius $1$ about the origin in $\Bbb{R}^{n+1}$, where $n\geqslant 2$. Choose $\epsilon>0$. Let $(\Sigma^n,\partial\Sigma^n)$ be a compact, locally convex hypersurface with boundary in $\Bbb{R}^{n+1}$, and suppose that the shape operator of $\Sigma$ is bounded below by $\epsilon\opId$ in the weak sense. Suppose, moreover, that $\Sigma$ passes through $0$. Let $H$ be a supporting hyperplane to $\Sigma$ at $0$. Let $d_H$ be the signed distance to $H$ in $\Bbb{R}^{n+1}$ such that $d_H$ is non-positive over $\Sigma$ in a neighbourhood of $0$. Let $(H_t)_{t\in\Bbb{R}}$ be the foliation of $\Bbb{R}^{n+1}$ by hyperplanes parallel to $H$. For all $t<0$, let $\Sigma_t$ be the connected component of $\Sigma$ lying above $H_t$ and containing $0$.
\proclaim{Proposition \nextprocno}
\noindent Choose $t<0$. Suppose that, for all $s\in]t,0[$, $\Sigma_t$ does not intersect either $\partial\Sigma$ or $\partial B_1(0)$. Then $\Sigma_t$ is embedded, and (together with $H_t$), $\Sigma_t$ bounds a convex set.
\endproclaim
\proclabel{PropBasicResult}
\proof We prove this by the method of moving planes, as in \cite{TrudWang}. In our case, strict convexity allows us to greatly simplify the argument, which we therefore include for the sake of clarity.
\medskip
\noindent We define $T\subseteq]t,0[$ such that $s\in T$ if and only if $\Sigma_s$ is embedded, and (together with $H_s$) bounds a convex set. For all $s$, we denote the convex set bounded by $\Sigma_s$ and $H_s$ by $K_s$. Since $\Sigma$ is strictly convex, all sufficiently small $s$ are elements of $T$, and $T$ is therefore non-empty. Let $t_0$ be the infimum of $T$. The result follows from the fact that all supporting hyperplanes to $\Sigma$ along $\Sigma_{t_0}\minter H_{t_0}$ are transverse to $H_{t_0}$, since, in this case, if $t<t_0$, then $\Sigma_{t_0}$ does not intersect either $\partial\Sigma$ or $\partial B_1(0)$, and $T$ may therefore be extended beyond $t_0$, which is absurd.
\medskip
\noindent To prove the assertion, we assume the contrary. Thus, choose $p\in\Sigma_{t_0}$ and suppose that the supporting hyperplane to $\Sigma_{t_0}$ at $p$ is tangent to $H_{t_0}$. Let $V_p$ be the unit normal vector to $H_{t_0}$ at $p$ which points outwards from $K_{t_0}$. $V_p=\pm\nabla d_H$. If $V_p=\nabla d_H$, then, since $K_{t_0}$ is convex, it is contained within $d_H^{-1}(]-\infty,t_0[)$, and therefore so is $0\in\Sigma_{t_0}$, which is absurd.
\medskip
\noindent Suppose now that $V_p=-\nabla d_H$. Let $(S^n,\partial S^n)$ be a compact $n$-dimensional topological manifold with boundary and $i:S\rightarrow\Bbb{R}^{n+1}$ a locally convex immersion such that $\Sigma=(S,i)$. Let $P\in S$ be the inverse image of $0$ in $S$ and, for all $t\in T$, let $S_t$ be the connected component of $(d_H\circ i)^{-1}(]t_0,+\infty[)$ containing $P$. $i$ restricts to a covering map from $S_{t_0}$ to $\partial K_{t_0}\setminus H_{t_0}$. The latter is however homeomorphic to a sphere with a point removed (i.e. a solid ball). Since this is simply connected, $i$ defines a homeomorphism from $S_{t_0}$ to $\partial K_{t_0}\setminus H_{t_0}$.
\medskip
\noindent Let $Q\in\partial S_{t_0}$ be such that $i(Q)=p$. Let $U_Q\subseteq S$ be a connected neighbourhood of $Q$ such that the restriction of $i$ to $U_Q$ is a homeomorphism onto an open subset of the boundary of a strictly convex set. Since $V_p=-\nabla d_H$ is the outward pointing supporting normal to $\Sigma$ at $Q$, by reducing $U_Q$ if necessary, we may assume that, throughout $U_Q\setminus\left\{Q\right\}$, $(d_H\circ i)>t_0$. Since $S$ has dimension at least $2$, $U_Q\setminus\left\{Q\right\}$ is connected, and therefore, by definition of $S_t$:
$$\matrix
&U\setminus\left\{Q\right\}\hfill&\subseteq S_t\hfill\cr
\Rightarrow\hfill&i(U_Q\setminus\left\{Q\right\})\hfill&\subseteq\partial K_{t_0}\hfill\cr
\Rightarrow\hfill&i(U_Q)\hfill&\subseteq\partial K_{t_0}.\hfill\cr
\endmatrix$$
\noindent By conservation of the domain, $i(U_Q)$ is an open subset of $\partial K_{t_0}$. $K_{t_0}$ is therefore strictly convex at $p$, and thus only meets $H_{t_0}$ at that point. Thus, by continuity, $i$ sends every point of $\partial S_{t_0}$ to $p$. Let $Q'$ be another point in $\partial S_{t_0}$. Let $U_{Q'}\subseteq S$ be a connected neighbourhood of $Q'$ such that the restriction of $i$ to $U_{Q'}$ is a homeomorphism onto an open subset of the boundary of a strictly convex set. Suppose, moreover, that:
$$
U_Q\minter U_{Q'} = \emptyset.
$$
\noindent Then, since $p\in i(U_Q)\minter i(U_{Q'})$:
$$
i(U_Q\setminus\left\{Q\right\})\minter i(U_{Q'}\setminus\left\{Q'\right\}) \neq\emptyset.
$$
\noindent However, the restriction of $i$ to $S_{t_0}$ is a homeomorphism, and it thus follows that $\partial S_{t_0}$ consists of only a single point, and $S$ is homeomorphic to $\partial K_{t_0}$, and is thus a topological sphere. In particular, $\partial S$ is trivial, which is absurd, and the claim follows.\qed
\medskip
\noindent Let $\gamma:]-T,T[\rightarrow\Bbb{R}^{n+1}$ be a unit speed curve whose geodesic curvature is bounded above by $\epsilon/(2+\epsilon)$. Suppose, moreover, that:
\medskip
\myitem{(i)} $\gamma$ is contained within $B_1(0)$;
\medskip
\myitem{(ii)} $\gamma$ lies in the exterior of $K_{t_0}$; and
\medskip
\myitem{(iii)} for all $t\in]-T,T[$, $(d_H\circ\gamma)(t)>t_0$.
\proclaim{Proposition \nextprocno}
\noindent If $d_K$ is the distance to $K_{t_0}$ in $\Bbb{R}^{n+1}$, then $(d_K+2/\epsilon)^2$ restricts to a convex function over $\gamma$.
\endproclaim
\proclabel{PropConvexityOfDistanceFunction}
\proof Choose $t\in]-T,T[$. Let $P\in\partial K_{t_0}$ be the closest point in $K_{t_0}$ to $\gamma(t)$. Since $\gamma$ lies above $H_{t_0}$, $P\in\Sigma_{t_0}$. Let $V_P$ be the outward pointing unit normal to $K_{t_0}$ at $P$ which is tangent to $\gamma$. We claim that the shape operator of $\Sigma$ with respect to $V_P$ is at least $(\epsilon/2)\opId$ in the weak sense. Indeed, if $V_P$ is a supporting normal to $\Sigma$, then the assertion follows by definition of $\Sigma$. Suppose, therefore, that $V_P$ is not a supporting normal to $\Sigma$. In particular, $P\in\partial\Sigma_{t_0}=:\Gamma$, and $V_P$ points upwards from $H_{t_0}$. Denote $V_1=-\nabla d_H$, and let $V_2$ be the unit supporting normal to $\Sigma$ lying in the plane spanned by $V_1$ and $V_P$. $V_1$ and $V_2$ define an obtuse, isosceles triangle whose angle at $p$ is $2\theta$, say. Since $V_P$ points upwards from $H_{t_0}$, the angle that $V_P$ makes with $V_1$ and $V_2$ respectively is $\theta+\varphi$ and $\theta-\varphi$, for some $0<\varphi<\theta$. If $A_1$ and $A_2$ are the shape operators of $\Gamma$ with respect to the normals $V_1$ and $V_2$ respectively, then:
$$
A_1 = 0,\qquad A_2\geqslant\epsilon\opId,
$$
\noindent in the weak sense. Thus, if $A_P$ is the shape operator of $\Gamma$ with respect to $P$, then:
$$
A_P \geqslant \frac{1}{2}(1 + \opSin(\varphi)/\opSin(\theta))A_2 > \frac{\epsilon}{2}\opId,
$$
\noindent in the weak sense. We may thus extend $\Gamma$ near $P$ to a convex hypersurface, $\Sigma_P$ such that:
\medskip
\myitem{(i)} $\Sigma_P$ lies outside $K_{t_0}$;
\medskip
\myitem{(ii)} $\Sigma_P$ meets $K_{t_0}$ at $P$;
\medskip
\myitem{(iii)} $V_P$ is the outward pointing unit normal to $\Sigma_P$ at $P$; and
\medskip
\myitem{(iv)} the shape operator of $\Sigma_P$ at $P$ is bounded below by $(\epsilon/2)\opId$ in the weak sense.
\medskip
\noindent The claim now follows. Define $r=(2/\epsilon)$. Let $Q\in\Bbb{R}^{n+1}$ be the unique point such that if $d_Q$ is the distance in $\Bbb{R}^{n+1}$ to $Q$, then:
\medskip
\myitem{(i)} $d_Q(P) = r$; and
\medskip
\myitem{(ii)} $(\nabla d_Q)(P) = V_P$.
\medskip
\noindent In a neighbourhood of $\gamma(t)$, $d_P^2\leqslant (d_K+(2/\epsilon))^2$. However:
$$\matrix
(\partial_t^2 d_P^2/2) \hfill&= (1/2)\opHess(d_P)(\partial_t\gamma,\partial_t\gamma) - d_P\langle\nabla d_P,\nabla_{\partial_t\gamma}\partial_t\gamma\rangle\hfill\cr
&\geqslant 1 - (1+r)\left|\nabla_{\partial_t\gamma}\partial_t\gamma\right|.\hfill\cr
\endmatrix$$
\noindent Moreover, by definition of the geodesic curvature:
$$
\left|\nabla_{\partial_t\gamma}\partial_t\gamma\right| \leqslant \frac{\epsilon}{2 + \epsilon} = \frac{1}{1+r}.
$$
\noindent The result follows.\qed
\proclaim{Proposition \nextprocno}
\noindent There exists $t_0<0$, which only depends on $\epsilon$, such that if, for all $s\in]t_0,0[$, $\Sigma_s$ does not intersect $\partial\Sigma$, then $\Sigma_{t_0}$ bounds a convex set $K_{t_0}$ (along with $H_{t_0}$) and $K_{t_0}$ is strictly contained within $B_1(0)$.
\endproclaim
\proclabel{CorBoundedAwayFromSphere}
\proof Without loss of generality, we may suppose that $\epsilon\leqslant 1$. Let $S$ be the sphere of radius $(4/\epsilon)$ which is tangent to $\Sigma$ at $0$ and locally contains $\Sigma$ in its interior. For all $t<0$, let $S_t$ be the connected component of $S$ lying above $H_t$ containing $0$. Let $t_0$ be such that $S_t$ is strictly contained inside $B_1(0)$.
\medskip
\noindent Let $T$ be the set of all $t<0$ such that, for all $s>t$, $\Sigma_t$ does not intersect either $\partial\Sigma$ or $\partial B_1(0)$. Since $\Sigma$ is strictly convex, for all $t$ small, $t\in T$, and $t$ is therefore non-empty. Let $t_1$ be the infimum of $T$ and suppose that $t_1>t_0$. By Proposition \procref{PropBasicResult}, $\Sigma_{t_1}$ bounds a convex set $K:=K_{t_1}$ (along with $H_{t_1}$). By the hypothesis on $\Sigma$, $\Sigma_{t_1}$ does not intersect $\partial\Sigma$, and thus $\Sigma_{t_1}$ intersects non-trivially with $\partial B_1(0)$, since otherwise $t_1$ would not be the infimum of $T$.
\medskip
\noindent Let $\gamma:\Bbb{R}\rightarrow S$ be a geodesic in $S$ such that $\gamma(0)=0$. Let $s_1$ be such that the connected component of $\gamma^{-1}(S_{t_1})$ containing $0$ coincides with $]-s_1,s_1[$. If $d$ is the distance in $\Bbb{R}^{n+1}$ to $K$, then $d\circ\gamma$ achieves a local minimum at $0$. Moreover, the geodesic curvature of $\gamma$, viewed as a curve in $\Bbb{R}^{n+1}$, is equal to $\epsilon/4<\epsilon/(\epsilon+2)$. Thus, by Proposition \procref{PropConvexityOfDistanceFunction}, the restriction of $d\circ\gamma$ to $]-s_1,s_1[$ is convex, and it follows that $d\circ\gamma$ is strictly positive over $]-s_1,s_1[\setminus\left\{0\right\}$. In particular, $\gamma(]-s_1,s_1[)$ lies outside $K$. Since $\gamma$ is arbitrary, the whole of $S_{t_1}$ lies outside $K$ and $K$ is therefore contained within the convex set bounded by $S_{t_1}$ and $H_{t_1}$. In particular, $K$ does not intersect $\partial B_1(0)$, which is absurd, and this completes the proof.\qed
\proclaim{Proposition \nextprocno}
\noindent For all $l\in[0,1]$, there exists $r>0$ which only depends on $l$ and $\epsilon$ such that, if $K\subseteq\Bbb{R}^{n+1}$ is a convex subset such that:
\medskip
\myitem{(i)} the shape operator of $\partial K$ is everywhere no less than $\epsilon\opId$ in the weak sense; and
\medskip
\myitem{(ii)} $K$ contains a geodesic segment of length $2l$,
\medskip
\noindent then $K$ contains a ball of radius $r$.
\endproclaim
\proclabel{PropNonTrivialInterior}
\proof Let $\gamma:[-l,l]\rightarrow K$ be a geodesic segment. Let $r_1=4/\epsilon$. Let $\eta:[-s,s]\rightarrow\Bbb{R}^{n+1}$ be a circular arc of radius at least $r_1$ and of angle no more than $\pi$ such that:
$$
\eta(-s) = \gamma(-l),\qquad \eta(s) = \gamma(l).
$$
\noindent Every point in $\eta$ lies at a distance no greater than $1$ from $\gamma$, and thus also from $K$. It follows by Proposition \procref{PropConvexityOfDistanceFunction} that $\eta$ lies entirely within $K$. Thus, if $\Omega\subseteq\Bbb{R}^{n+1}$ is the locus of all points traced out by such arcs, then:
$$
\Omega\subseteq K.
$$
\noindent However, if $r>0$ is such that:
$$
(r_1 - r)^2 + l^2 = r_1^2.
$$
\noindent Then:
$$
B_r(\gamma(0)) \subseteq \Omega.
$$
\noindent The result follows.\qed
\goodbreak
\newhead{Strictly Convex Hypersurfaces of General Manifolds}
\noindent Let $M^{n+1}$ be an $(n+1)$-dimensional Riemannian manifold. Let $K\subseteq M$ be a compact subset. Choose $\epsilon>0$, and let $(\Sigma,\partial\Sigma)$ be a compact, locally convex immersed hypersurface in $M$ such that:
\medskip
\myitem{(i)} $\Sigma$ is contained in $K$;
\medskip
\myitem{(ii)} $\Gamma:=\partial\Sigma$ is smooth; and
\medskip
\myitem{(iii)} the shape operator of $\Sigma$ is everywhere greater than $\epsilon\opId$ in the weak sense.
\proclaim{Lemma \nextprocno}
\noindent There exists $r_1>r_2>0$ (which only depend on $\epsilon$, $M$, $K$ and $\Gamma$) such that, for all $p\in\Sigma$:
\medskip
\myitem{(i)} the connected component of $\Sigma\minter B_{r_1}(p)$ containing $p$ is embedded and lies on the boundary of a convex set; and
\medskip
\myitem{(ii)} moreover, this convex set contains an open ball of radius $r_2$.
\endproclaim
\proclabel{LemmaPrecompactnessOfStrictlyConvexImmersions}
\noindent Choose $\epsilon>0$ and let $(\Sigma',\partial\Sigma')$ be a smooth, strictly convex, immersed hypersurface in $\Bbb{R}^{n+1}$ such that:
\medskip
\myitem{(i)} the shape operator of $\Sigma'$ is greater than $\epsilon$ in the weak sense; and
\medskip
\myitem{(ii)} $\partial\Sigma'$ consists of two connected components, which we denote by $\partial\Sigma_1'$ and $\partial\Sigma_2'$ respectively.
\medskip
\noindent If $\partial\Sigma_1'=\Gamma$, then we refer to $(\Sigma',\partial\Sigma')$ as a {\bf thickening} of $\Gamma$. Let $\msf{N}$ be the outward pointing unit normal vector field over $\Sigma'$. Choose $P\in\partial\Sigma_2'$. Let $\Cal{N}_P$ denote the circle of unit vectors normal to $\partial\Sigma_2'$ at $P$. For any vector, $V$ in $\Cal{N}_P$, let $H_V$ be the (oriented) hyperplane in $T_PM$ normal to $V$ at $P$. We identify $H_V$ with its image under the exponential map of $M$, and we define $\Sigma_V'$ to be the connected component of $\Sigma'$ lying above $H_V$ and containing $P$. We have the following elementary result:
\goodbreak
\proclaim{Proposition \nextprocno}
\myitem{(i)} There exists $\theta_{\Sigma'}>0$, which only depends on $\Sigma'$, such that, if the angle between $V$ and $\msf{N}(P)$ is less than $\theta_{\Sigma'}$, then $\Sigma'_V$ does not intersect $\partial\Sigma'_1$; and
\medskip
\myitem{(ii)} there exists $h_{\Sigma'}>0$, which also only depends on $\Sigma'$ such that, if $V$ makes an angle of at least $\theta_{\Sigma'}$ with $\msf{N}(P)$, then the furthest point in $\Sigma'_V$ from $H_V$ lies at a distance no less than $h_{\Sigma'}$ from $H_V$.
\endproclaim
\proclabel{PropBoundsObtainedFromTheBoundary}
\proof If $V$ coincides with $\msf{N}(P)$, then $\Sigma'_V=\left\{P\right\}$. The first assertion follows by continuity and compactness of $\partial\Sigma_2'$. The second assertion now follows by compactness of $\partial\Sigma_2'$, and this completes the proof.\qed
\medskip
{\bf\noindent Proof of Lemma \procref{LemmaPrecompactnessOfStrictlyConvexImmersions}:\ }Let $S$ be a thickening of $\Gamma$ such that:
\medskip
\myitem{(i)} $S$ is homeomorphic to $\Gamma\times[0,1]$;
\medskip
\myitem{(ii)} if $\pi:S\rightarrow[0,1]$ is the canonical projection, then $\pi^{-1}(\left\{0\right\})$ coincides with $\Gamma$;
\medskip
\myitem{(iii)} the shape operator of $S$ is bounded below by $\epsilon/2$; and
\medskip
\myitem{(iv)} $S\munion\Sigma$ is a locally convex immersed hypersurface.
\medskip
\noindent $S$ may be chosen independant of $\Sigma$ (despite condition $(iv)$). Indeed, if $S$ is chosen such that:
\medskip
\myitem{(i)} its outward pointing normal along $\Gamma$ points in the same direction as that of $\Sigma$; and
\medskip
\myitem{(ii)} the lowest eigenvalue of the shape operator of $S$ along $\Gamma$ is always in $[\epsilon/4,\epsilon/2]$,
\medskip
\noindent then $S$ satisfies condition $(iv)$ (see Section $4$ of \cite{SmiPPH} for details).
\medskip
\noindent Let $g$ be the Riemannian metric over $M$. For $q\in M$, identify $T_q M$ with $\Bbb{R}^{n+1}$, and let $\opExp_q:\Bbb{R}^{n+1}\rightarrow M$ be the exponential map. Let $\nabla^0$ be the Euclidean covariant derivative over $\Bbb{R}^{n+1}$, let $\nabla^q$ be the pull-back through $\opExp_q$ of the Levi-Civita covariant derivative of $g$, and let $\Omega^q=\nabla^q-\nabla^0$ be the connexion $2$-form of $\nabla^q$ with respect to $\nabla^0$. There exists $r_1>0$, which only depends on $K$ such that, for all $q\in K$, $\|\Omega^q\|<\epsilon/4$ over the ball of radius $r_1$ about $q$. If we replace $g$ with $r_1^{-1}g$, then the shape operator of $S\munion\Sigma$ is bounded below by $(\epsilon r_1)/2$ in the weak sense, and $\|\Omega^q\|<(\epsilon r_1)/4$ over the ball of radius $1$ about $q$, for all $q\in K$.
\medskip
\noindent Denote $\epsilon_1 = (\epsilon r_1)/4$. Choose $p\in\Sigma$, and let $(\tilde{\Sigma},\partial\tilde{\Sigma})$ be the connected component of $S\munion\Sigma$ lying in $B_1(p)$ containing $p$. We identify $(\tilde{\Sigma},\partial\tilde{\Sigma})$ with its pull back through $\opExp_q$ in $\Bbb{R}^{n+1}$, and observe that the shape operator of $\tilde{\Sigma}$ with respect to the rescaled Euclidean metric over $\Bbb{R}^{n+1}$ is everywhere bounded below by $\epsilon_1$ in the weak sense.
\medskip
\noindent Let $H$ be a supporting hyperplane to $\tilde{\Sigma}$ at $p$. Let $d$ be the (signed) distance in $\Bbb{R}^{n+1}$ to $H$ such that, near $p$, $d$ is non-positive over $\tilde{\Sigma}$. For all $t\in\Bbb{R}$, let $H_t$ be the hyperplane parallel to $H$ lying at distance $t$ from $H$. For all $t<0$, let $\tilde{\Sigma}_t$ be the connected component of $\tilde{\Sigma}$ lying above $H_t$ and containing $p$. Let $t_0<0$ be as in Proposition \procref{CorBoundedAwayFromSphere}. Let $\theta_S$ and $h_S$ be as in Proposition \procref{PropBoundsObtainedFromTheBoundary}, now chosen such that the proposition remains valid with respect to the rescaled Euclidean metric over $\Bbb{R}^{n+1}$ (as opposed to $g$). Suppose that $\left|t_0\right|<h_S$.
\medskip
\noindent Let $T$ be the set of all $t\in]t_0,0[$ such that for all $s\in]t,0[$, $\tilde{\Sigma}_t$ does not intersect either $\partial\tilde{\Sigma}$ or $\partial B_1(0)$. For all $t$ sufficiently close to $0$, $t\in T$. In particular, $T$ is non-empty. Let $t_1$ be the infimum of $T$. By Proposition \procref{PropBasicResult}, $\Sigma_{t_1}$ (along with $H_{t_1}$) bounds a convex set, $K$, say.
\medskip
\noindent Suppose now that $t_1>t_0$. By Proposition \procref{PropBasicResult}, since $t_1$ is the infimum of $T$, $\tilde{\Sigma}_{t_1}$ meets either $\partial\tilde{\Sigma}$ or $\partial B_1(p)$, since, otherwise, it could be extended further. By Proposition \procref{CorBoundedAwayFromSphere}, $\tilde{\Sigma}_{t_1}$ does not intersect $\partial B_1(0)$. It follows that $\tilde{\Sigma}_{t_1}$ intersects $\partial\tilde{\Sigma}_{t_1}$ at some point, $Q$, say. $H_{t_1}$ is tangent to $\partial\tilde{\Sigma}_{t_1}$ at $Q$. Since $\tilde{\Sigma}_{t_1}$ intersects non-trivially with $\partial\Sigma$, it follows by assertion $(i)$ of Proposition \procref{PropBoundsObtainedFromTheBoundary} that $H_{t_1}$ makes an angle of at least $\theta_S$ with the outward pointing normal to $H_S$ at $Q$. Since the connected component of $S$ lying above $H_{t_1}$ is contained in $\tilde{\Sigma}_{t_1}$, it follows by assertion $(ii)$ of Proposition \procref{PropBoundsObtainedFromTheBoundary} that the supremum of the distance to $H_{t_0}$ over $\tilde{\Sigma}_{t_1}$ is at least $h_S$. However, by convexity, $\tilde{\Sigma}_{t_1}$ lies below $H_0$, and so:
$$
h_S \leqslant \left|t_1\right| < \left|t_0\right| < h_S.
$$
\noindent This is absurd, and it follows that $t_1=t_0$. Assertion $(i)$ now follows for $r_1<\left|t_0\right|$. Assertion $(ii)$ follows by Proposition \procref{PropNonTrivialInterior}, since $\Sigma_{t_0}\minter\partial B_{r_1}(p)\neq\emptyset$. This completes the proof.\qed
\goodbreak
\newhead{Compactness}
\noindent Let $M^{n+1}$ be an $(n+1)$-dimensional Riemannian manifold. Choose $\epsilon>0$ and let $(\Sigma_n,\partial\Sigma_n)_\ninn$ be a sequence of compact, locally convex immersed hypersurfaces in $M$ such that, for all $n$:
\medskip
\myitem{(i)} $\Gamma_n:=\partial\Sigma_n$ is $C^\infty$; and
\medskip
\myitem{(ii)} the shape operator of $\Sigma_n$ is greater than $\epsilon\opId$ in the weak sense.
\medskip
\noindent Let $\Gamma_0$ be a smooth, compact, codimension $2$ immersed submanifold of $M$ and suppose that $(\Gamma_n)_\ninn$ converges to $\Gamma_0$ in the $C^\infty$ sense. In other words, for sufficiently large $n$, $\Gamma_n$ is a normal graph over $\Gamma_0$ and the (unique) function of which $\Gamma_n$ is a graph tends to $0$ in the $C^\infty$ sense.
\medskip
\noindent For all $n$, let $\opVol(\Sigma_n)$ be the volume of $\Sigma_n$ with respect to the intrinsic measure.
\proclaim{Lemma \nextprocno}
\noindent Suppose that there exists a compact subset $K\subseteq M$ and a real number $B>0$ such that, for all $n\in\Bbb{N}$:
$$
\Sigma_n\subseteq K,\qquad \opVol(\Sigma_n)\leqslant B.
$$
\noindent Then there exists a locally convex immersed hypersurface $\Sigma_0$ in $M$ such that:
\medskip
\myitem{(i)} $\Gamma_0=\partial\Sigma_0$;
\medskip
\myitem{(ii)} the shape operator of $\Sigma_0$ is greater than $\epsilon\opId$ in the weak sense; and
\medskip
\myitem{(iii)} $(\Sigma_n)_\ninn$ subconverges to $\Sigma_0$.
\endproclaim
\proclabel{LemmaCompactness}
\proof For all $n$, let $\Sigma^c_n$ be a thickening of $\Gamma_n$ such that:
\medskip
\myitem{(i)} $\Sigma^c_n$ is homeomorphic to $\Gamma_n\times[0,1]$;
\medskip
\myitem{(ii)} if $\pi_n:\Sigma^c_n\rightarrow[0,1]$ is the canonical projection, then $\pi^{-1}(\left\{0\right\})$ coincides with $\Gamma_n$;
\medskip
\myitem{(iii)} the shape operator of $\Sigma^c_n$ is bounded below by $\epsilon/2$; and
\medskip
\myitem{(iv)} $\hat{\Sigma}_n:=\Sigma^c_n\munion\Sigma_n$ is a locally convex immersed hypersurface.
\medskip
\noindent Suppose, moreover, that $(\Sigma^c_n)_\ninn$ converges to $\Sigma^c_0$ in the $C^\infty$ sense. As in the proof of Lemma \procref{LemmaPrecompactnessOfStrictlyConvexImmersions}, the $(\Sigma_n^c)_\ninn$ may be chosen independant of $(\Sigma_n)_\ninn$ (despite condition $(iv)$). For all $n$, let $S_n\subseteq\hat{S}_n$ be abstract manifolds and $i_n:\hat{S}_n\rightarrow M$ be a locally convex immersion such that:
$$
\Sigma_n = (i_n,S_n),\qquad \hat{\Sigma}_n = (i_n,\hat{S}_n).
$$
\noindent For all $n$, we furnish $\hat{S}_n$ with the metric and the measure induced by $i_n$. For all $n$, for all $p\in S_n$ and for all $\epsilon>0$, let $B_{n,\epsilon}(p)$ be the intrinsic ball of radius $\epsilon$ about $p$ in $\hat{\Sigma}_n$.
\medskip
\noindent Let $r_1>r_2>0$ be as in Lemma \procref{LemmaPrecompactnessOfStrictlyConvexImmersions}.  Choose $\epsilon\in]0,r_1[$. For all $n$, let $p_n$ be a point in $\Sigma_n$. By Lemma \procref{LemmaPrecompactnessOfStrictlyConvexImmersions}, for all $n$, there exists a convex set $K_n\subseteq M$ such that:
\medskip
\myitem{(i)} $K_n$ contains an open ball of radius $r_2$; and
\medskip
\myitem{(ii)} the restriction of $i_n$ to $B_{n,r_1}(p_n)$ is a homeomorphism onto the open ball of radius $r_1$ about $i_n(p_n)$ in $\partial K_n$.
\medskip
\noindent By compactness of the family of compact, convex sets, there exists a compact, convex set $K_0$ towards which $(K_n)_\ninn$ subconverges. Moreover, by $(i)$, $K_0$ contains an open ball of radius $r_2$, and therefore has non-trivial interior. $i_n(B_{n,\epsilon}(p_n))$ subconverges to an open ball of radius $\epsilon$ about some point in $\partial K_0$. In particular $(\opVol(B_{n,\epsilon}(p_n)))_\ninn$ tends towards the volume of this ball, which is non-zero. Since $(p_n)_\ninn$ was arbitrary, it follows that, for all $\epsilon>0$, there exists $v_\epsilon$ such that, for all $n$ and for all $p\in S_n$:
$$
\opVol(B_{n,\epsilon}(p))\geqslant v_\epsilon.
$$
\noindent By increasing $B$ if necessary, we may suppose that, for all $n$, $\opVol(\hat{S}_n)\leqslant B$. Thus, if $N_\epsilon\geqslant B/v_\epsilon$, then, for all $n$, the maximum number of disjoint balls in $\hat{S}_n$ of radius $\epsilon$ with centres in $S_n$ is at most $N_\epsilon$. Thus, the minimum number of balls in $\hat{S}_n$ of radius $2\epsilon$ required to cover $S_n$ is at most $N_\epsilon$. It follows by the fundamental compactness theorem of metric spaces that there exists a metric space $S_0$ towards which $(S_n)_\ninn$ subconverges in the Gromov/Hausdorff sense (see \cite{Petersen}).
\medskip
\noindent Choose $p\in S_0$. For all $n$, let $p_n$ be a point in $S_n$ and suppose that $(p_n)_\ninn$ converges to $p_0$. Constructing $(K_n)_\ninn$ and $K_0$ as before, we obtain an isometry, $i_0$, from $B_{0,\epsilon}(p_0)$ to a ball of radius $\epsilon$ about a point in $\partial K_0$. Moreover, the sequence of restrictions of $(i_n)_\ninn$ to the balls of radius $\epsilon$ about the $(p_n)_\ninn$ converges locally uniformly to $i_0$. Using a diagonal argument, we obtain a $C^{0,1}$ mapping $i_0$ from the whole of $S_0$ into $M$ which is a locally convex immersion. $\Sigma_0=(S_0,i_0)$ is the desired hypersurface, and this completes the proof.\qed
\goodbreak
\newhead{The Dirichlet Problem}
\noindent Let $M^{n+1}$ be an $(n+1)$-dimensional Riemannian manifold of non-positive curvature. Choose $\theta\in[(n-1)\pi/2,n\pi/2[$. Let $H\subseteq M$ be a smooth, locally convex hypersurface. Let $\Omega\subset H$ be a bounded, open subset of $H$. Let $\hat{\Sigma}\subseteq M$ be a convex hypersurface such that $\partial\hat{\Sigma}=\partial\Omega=:\Gamma$. Let $\phi:M\rightarrow]0,\infty[$ be a smooth, positive function such that, for all $p\in H$:
$$
R_\theta(H)(p) < \phi(p),
$$
\noindent and, for all $p\in\hat{\Sigma}$:
$$
R_\theta(\hat{\Sigma})(p) > \phi(p),
$$
\noindent in the weak sense.
\medskip
\noindent Theorem $3.22$ of \cite{SmiNLD} may be adapted to yield:
\proclaim{Proposition \nextprocno}
\noindent Suppose that $\hat{\Sigma}$ is a graph over $\Omega$ and that $\Gamma$ is strictly convex as a subset of $M$ with respect to the outward pointing normal to $\Gamma$ in $\hat{\Sigma}$. If $\theta>(n-1)\pi/2$, then there exists an immersed hypersurface $\Sigma\subseteq M$ such that:
\medskip
\myitem{(i)} $\Sigma$ is $C^0$ and $C^\infty$ in its interior;
\medskip
\myitem{(ii)} $\partial\Sigma = \Gamma$;
\medskip
\myitem{(iii)} $\Sigma$ is a graph over $\Omega$ lying below $\hat{\Sigma}$; and
\medskip
\myitem{(iv)} for all $p\in\Sigma$, $\hat{R}_\theta(\Sigma)(p) = \phi(p)$.
\medskip
\noindent Moreover, the same result holds for $\theta=(n-1)\pi/2$ provided that, in addition, the shape operator of $\hat{\Sigma}$ is everywhere bounded below by $\epsilon\opId$ in the weak sense, for some $\epsilon>0$.
\endproclaim
\proclabel{PropBasicExistenceResult}
\remark In fact, we don't need the result for $\theta=(n-1)\pi/2$.
\medskip
\remark The hypotheses of this proposition are satisfied when the norm of the second fundamental form of $H$ is small with respect to that of $\hat{\Sigma}$ and the normal of $\hat{\Sigma}$ is sufficiently bounded away from $TH$ along $\Gamma$. Explicitely, if the second fundamental form of $\hat{\Sigma}$ is greater than $\epsilon\opId$ in the weak sense, if the norm of the second fundamental form of $H$ is bounded above by $\delta$, and if the angle between the normal to $\hat{\Sigma}$ and $TH$ is bounded below by $\theta$ along $\Gamma$, then the hypotheses are satisfied provided that:
$$
\epsilon\opSin(\theta) - \delta > 0.
$$
{\noindent\bf Sketch of proof:\ }We use the continuity method, which reduces to showing openness and closedness of a certain interval, from which existence is deduced by connectedness. In \cite{SmiNLD}, openness is shown by proving that the linearisation of the special-Lagrangian curvature operator is always invertible, which follows from the stronger hypotheses used there by Lemma $7.5$ of \cite{SmiFCS}. In our case, the hypotheses on the curvature of the ambient manifold and the special Lagrangian curvature of the immersed hypersurface are weaker, but openness may still be obtained (generically), using the differential topological argument discussed (in the case of Gaussian curvature) in Section $11$ of \cite{SmiPPH}. This approach requires a strengthening of the compactness result (used to prove closedness) to treat the case of hypersurfaces of prescribed (non-constant) curvature. This strengthening is required anyway to obtain the existence result in the generality given here, and is obtained by a fairly trivial adaptation of the reasoning presented in Sections $3.1$ to $3.6$ of \cite{SmiNLD} used there to prove compactness in the case of hypersurfaces of constant special Lagrangian curvature. This completes the proof.\qed
\goodbreak
\newhead{The Plateau Problem}
\noindent Let $M^{n+1}$ be a Hadamard manifold. Choose $\theta\in[(n-1)\pi/2,n\pi/2[$. Let $(\hat{N}^n,\partial\hat{N}^n)$ be a locally convex immersed hypersurface in $M$. Let $\phi:M\rightarrow]0,\infty[$ be a strictly positive, smooth function such that, for all $p\in\hat{N}$:
$$
R_\theta(N)(p) \geqslant \phi(p),
$$
\noindent in the weak sense.
\medskip
{\bf\noindent Proof of Theorem \procref{ThmExistence}:\ }We first consider the case where $\theta>(n-1)\pi/2$. By Lemma $2.2$ of \cite{SmiNLD}, there therefore exists $k>0$ such that, if $N'$ is a smooth, locally convex immersed hypersurface, and if $A$ is the shape operator of $N'$, then:
$$
A \geqslant kR_\theta(N').
$$
\noindent Now let $(N,\partial N)$ be an $n$-dimensional, locally convex immersed hypersurface in $M$ and suppose that, for all $p\in N$:
$$
R_\theta(N)(p) \geqslant \phi(p),
$$
\noindent in the weak sense. Let $p\in N$ be an interior point. Let $\msf{N}_p$ be a supporting normal to $N$ at $p$. By Lemma $4.7$ of \cite{SmiNLD}, we may suppose that there exists $\eta>0$ such that, if $\msf{N}'_p$ is any other supporting normal to $N$ at $p$, then:
$$
\langle \msf{N}_p,\msf{N}'_p\rangle \geqslant \eta.
$$
\noindent Choose $\epsilon>0$ such that $\epsilon< kR_\theta(N)$, and let $H$ be a smooth, strictly convex hypersurface passing through $p$ such that:
\medskip
\myitem{(i)} the outward pointing unit normal to $H$ at $p$ coincides with $\msf{N}_p$; and
\medskip
\myitem{(ii)} the norm of the shape operator of $H$ is always less than $\epsilon$.
\medskip
\noindent We extend $H$ to a foliation $(H_t)_{t\in]-\delta,\delta[}$ of a neighbourhood of $p$ such that norm of the shape operator of each leaf of the foliation is always less than $\epsilon$. For all sufficiently small $t<0$, let $N_t$ be the connected component of $N$ lying above $H_t$ and containing $p$. For $\epsilon$ sufficiently small, there exists $t_0<0$, and, for all $t\in]t_0,0[$, an open subset $\Omega_t\subseteq H_t$ such that $N_t$ is a graph over $\Omega_t$. Moreover, by reducing $\epsilon$ further if necessary, for all $t\in]t_0,0[$, $H_t$ and $N_t$ satisfy the hypotheses of Proposition \procref{PropBasicExistenceResult}. There thus exists, for all such $t$, a smooth, locally  convex immersed hypersurface $N'_t$ lying between $H_t$ and $N_t$ such that, for all $p\in N'_t$:
$$
R_\theta(N'_t)(p) = \phi(p).
$$
\noindent Let $N'$ be the locally convex immersed hypersurface obtained by replacing $N_{t_0}$ with $N'_{t_0}$. For all $p\in N'$ lying in the complement of $\partial N'_{t_0}$:
$$
R_\theta(N')(p) \geqslant \phi(p),
$$
\noindent in the weak sense. By deforming slightly, we may suppose that $N'_{t_0}$ is smooth up to the boundary. By Lemma $2.4$ of \cite{SmiNLD}, since $R_\theta$ is a concave function, it follows that $R_\theta(N')(p)\geqslant\phi(p)$ in the weak sense also for all $p\in\partial N'_{t_0}$, and therefore over the whole of $N'$. We observe finally that $N'$ is trivially bounded by $N$, and we call $N'$ a {\bf Perron Regularisation} of height $t_0$ of $N$ about $p$.
\medskip
\noindent Let $\Cal{F}$ denote the family of all locally strictly convex, immersed hypersurfaces that may be obtained from $N$ by a finite number of iterations of the Perron Regularisation process. For all $N\in\Cal{F}$, let $\opVol(N)$ be the volume of $N$. $\opVol$ defines a continuous function from $\Cal{F}$ into $[0,\infty[$. Observe that, by convexity, and since the ambient manifold is non-positively curved, Perron Regularisation always reduces volume. Define $V_0$ by:
$$
V_0 = \minf\left\{ \opVol(N)\text{ s.t. }N\in\Cal{F}\right\}.
$$
\noindent Let $(N_n)_\ninn\in\Cal{F}$ be such that:
$$
\opVol(N_n)_\ninn\rightarrow V_0.
$$
\noindent By Lemma \procref{LemmaCompactness}, $(N_n)_\ninn$ subconverges to a locally strictly convex $C^{0,1}$ immersed hypersurface, $N_0$ such that:
$$
\opVol(N_0) = V_0.
$$
\noindent We aim to show that $N_0$ is smooth away from the boundary, and at every point $p\in N_0$:
$$
R_\theta(N_0)(p) = \phi(p).
$$
\noindent Choose $p_0\in N_0$. For all $n$, choose $p_n\in N_n$ such that $(p_n)_\ninn$ converges to $p_0$. Since $N_0$ is locally strictly convex, and since $(N_n)_\ninn$ converges to $N_0$, there exists $t_0<0$, and, for all $n\in\Bbb{N}$, a Perron Regularisation $N_n'$ of height $t_0$ of $N_n$ about $p_n$. For all $n$, $N_n'$ is bounded by $N_n$ and $\opVol(N_n')\leqslant\opVol(N_n)$. By Lemma \procref{LemmaCompactness}, $(N_n')_\ninn$ subconverges to a locally convex immersed hypersurface, $N'_0$, say.  By continuity, $\opVol(N_0')\leqslant\opVol(N_0)$. Moreover, by Lemma \procref{LemmaClosednessOfGraphProperty}, $N_0'$ is bounded by $N_0$. However, by definition of $V_0$:
$$
\opVol(N_0')\geqslant V_0 = \opVol(N_0).
$$
\noindent Consequently $\opVol(N_0')=\opVol(N_0)$, and since $N_0'$ is bounded by $N_0$, the two hypersurfaces must coincide. For all $n$, choose $p'_n\in N'_n$ such that $(p'_n)_\ninn$ converges to $p_0$. For all $r>0$ and for all $n\in\Bbb{N}\munion\left\{0\right\}$, let $B_r'(p'_n)$ be the ball of radius $r$ (with respect to the intrinsic distance) about $p'_n$ in $N'_n$. Choose $r_0<\left|t_0\right|$. For all $n$, $N'_n$ is smooth over the ball of radius $r_0$ about $p'_n$, and for all $p\in B_{r_0}'(p_n')$:
$$
R_\theta(N'_n)(p) = \phi(p).
$$
\noindent Taking limits, it follows by Theorem $1.4$ of \cite{SmiSLC} that $N_0'=N_0$ is also smooth over the ball of radius $r$ about $p_0$, and that, for all $p\in B_r(p_0)$:
$$
R_\theta(N_0)(p) = \phi(p).
$$
\noindent Since $p_0\in N_0$ is arbitrary, the result now follows for the case where $\theta>(n-1)\pi/2$.
\medskip
\noindent Now suppose that $\theta=(n-1)\pi/2$. Let $(\delta_n)_\ninn>0$ be a sequence converging to $0$, let $(\theta_n)_\ninn>(n-1)\pi/n$ be a sequence converging to $\theta$ and, for all $n\in\Bbb{N}$, let $N_n$ a locally strictly convex hypersurface with smooth interior such that $\partial N_n = \Gamma$ and, for every interior point $p\in N_n$:
$$
R_{\theta_n}(N_n)(p) = \phi(p) - \delta_n.
$$
\noindent This is possible, since the shape operator of $N$ is everywhere bounded below by $\epsilon\opId$ in the weak sense, for some $\epsilon>0$.
\medskip
\noindent For $\epsilon>0$ and $n\in\Bbb{N}$, let $N_{\epsilon,n}$ be the set of all points in $N_n$ whose (intrinsic) distance from $\partial N_n$ is at most $\epsilon$. By Lemma $5.1$ of $\cite{SmiPPH}$, there exists $\epsilon>0$ such that $N_{\epsilon,n}$ subconverges to a locally strictly convex hypersurface, $N_{\epsilon,0}$, say. Choose $r>0$ and for all $n\in\Bbb{N}$, let $p_n\in N_n$ be a point whose (intrinsic) distance from $\partial N_n$ is at least $r$ and let $B_{n,r}(p_n)$ be the (intrinsic) ball of radius $r$ about $p_n$. By Theorem $1.4$ of $\cite{SmiSLC}$, $(B_{n,r}(p_n))_\ninn$ converges either to a smooth immersed hypersurface or to a geodesic segment of length $r$.
\medskip
\noindent For all $n$, since $N_n$ is bounded by $N$, $N$ is (more or less) a graph over $N_n$, and there exists a canonical projection $\pi_n:N\rightarrow N_n$ (see Section $3$ of \cite{SmiPPH} for details). For all $n$, $\pi_n$ is $1$-Lipschitz and coincides with the identity along $\partial N=\partial N_n$. Thus, by the Arzela-Ascoli Theorem, we may assume that there exists $\pi_0:N\rightarrow M$ towards which $(\pi_n)_\ninn$ converges. The set of all $p\in N$ where $(B_{n,r}(\pi_n(p)))_\ninn$ converges to a smooth immersed hypersurface is open, likewise so is the set where it converges towards a geodesic segment. Since $N_\epsilon$ is contained in the former of the two, it follows that this subset is non-trivial, and thus, by connectedness, coincides with the whole of $N$. Thus every such limit is a smooth immersed hypersurface, and we deduce, as before, that there exists a convex immersed hypersurface $N_0$ towards which $(N_n)_\ninn$ subconverges. Moreover, by Theorem $1.4$ of $\cite{SmiSLC}$, $N_0$ is smooth away from the boundary, and for every interior point $p\in N_0$:
$$
R_\theta(N_0)(p) = \phi(p).
$$
\noindent This completes the proof.\qed
\goodbreak
\newhead{Bibliography}
{\leftskip = 5ex \parindent = -5ex
\leavevmode\hbox to 4ex{\hfil \cite{AndBarbBegZegh}}\hskip 1ex{Andersson L., Barbot T., B\'eguin F., Zeghib A., Cosmological time versus CMC time in spacetimes of constant curvature, arXiv:math/0701452}
\medskip
\leavevmode\hbox to 4ex{\hfil \cite{Ghomi}}\hskip 1ex{Alexander S., Ghomi M., Wong J., Topology of Riemannian submanifolds with prescribed boundary, {\sl to appear in Duke Math. J.}}
\medskip
\leavevmode\hbox to 4ex{\hfil \cite{GuanSpruck}}\hskip 1ex{Guan B., Spruck J., The existence of hypersurfaces of constant Gauss curvature with prescribed boundary, {\sl J. Differential Geom.} {\bf 62} (2002), no. 2, 259--287}
\medskip
\leavevmode\hbox to 4ex{\hfil \cite{LabA}}\hskip 1ex{Labourie F., Un lemme de Morse pour les surfaces convexes (French), {\sl Invent. Math.} {\bf 141} (2000), no. 2, 239--297}
\medskip
\leavevmode\hbox to 4ex{\hfil \cite{Petersen}}\hskip 1ex{Petersen P., {\sl Riemannian geometry}, Graduate Texts in Mathematics, {\bf 171}, Springer, New York, (2006)}
\medskip
\leavevmode\hbox to 4ex{\hfil \cite{Ros}}\hskip 1ex{Rosenberg H., Hypersurfaces of constant curvature in space forms, {\sl Bull. Sci. Math.} {\bf 117} (1993), no. 2, 211--239}
\medskip
\leavevmode\hbox to 4ex{\hfil \cite{SmiSLC}}\hskip 1ex{Smith G., Special Lagrangian curvature, arXiv:math/0506230}
\medskip
\leavevmode\hbox to 4ex{\hfil \cite{SmiFCS}}\hskip 1ex{Smith G., Moduli of Flat Conformal Structures of Hyperbolic Type, arXiv:0804.0744}
\medskip
\leavevmode\hbox to 4ex{\hfil \cite{SmiNLD}}\hskip 1ex{Smith G., The Non-Linear Dirichlet Problem in Hadamard Manifolds,\break arXiv:0908.3590}
\medskip
\leavevmode\hbox to 4ex{\hfil \cite{SmiCGC}}\hskip 1ex{Smith G., Constant Gaussian Curvature Hypersurfaces in Hadamard Manifolds,\break arXiv:0912.0248}
\medskip
\leavevmode\hbox to 4ex{\hfil \cite{SmiPPH}}\hskip 1ex{Smith G., The Plateau Problem in Hadamard Manifolds, arXiv:1002.2982}
\medskip
\leavevmode\hbox to 4ex{\hfil \cite{TrudWang}}\hskip 1ex{Trudinger N. S., Wang X., On locally locally convex hypersurfaces with boundary, {\sl J. Reine Angew. Math.} {\bf 551} (2002), 11--32}
\par}
\enddocument